\documentclass[11pt]{amsproc}
\usepackage[margin=1in]{geometry}
\usepackage{setspace,fullpage}
\usepackage{graphicx}
\usepackage[nice]{nicefrac}
\usepackage{amssymb}
\usepackage{multirow}
\usepackage{array}

\usepackage{amsmath,amsthm,amscd,amssymb, mathrsfs}
\usepackage{latexsym}

\numberwithin{equation}{section}

\theoremstyle{plain}
\newtheorem{theorem}{Theorem}[section]
\newtheorem{lemma}[theorem]{Lemma}

\newtheorem{proposition}[theorem]{Proposition}
\newtheorem{conjecture}[theorem]{Conjecture}

\theoremstyle{definition}
\newtheorem{definition}[theorem]{Definition}

\theoremstyle{remark}

\newtheorem{case[theorem]}{Case}

\title[\parbox{14cm}{\centering{ $\mathcal{S}-operators\hspace{1in}}} \quad]{Connections between $\mathcal{S}$-operators and restriction estimates for spheres over finite fields}
\author{Hunseok Kang and Doowon Koh }
 
\address{College of Engineering and Technology\\
American University of the Middle East\\
Kuwait}
\email{hunseok.kang@aum.edu.kw}

\address{Department of Mathematics\\
Chungbuk National University \\
Cheongju, Chungbuk 28644 Korea}
\email{koh131@chungbuk.ac.kr}

\thanks{
 Doowon Koh was supported by the National Research Foundation of Korea
(NRF) grant funded by the Korea government (MSIT) (NO. RS-2023-00249597).}

\begin{document} 

\begin{abstract} In this paper, we introduce a new operator, $\mathcal{S}$, which is closely related to the restriction problem for spheres in $\mathbb{F}_q^d$, the $d$-dimensional vector space over the finite field $\mathbb{F}_q$ with $q$ elements. The $\mathcal{S}$ operator is considered as a specific operator that maps functions on $\mathbb{F}_q^d$ to functions on $\mathbb{F}_q^{d+1}$.
We explore a relationship between the boundedness of the $\mathcal{S}$ operator and the restriction estimate for spheres in $\mathbb{F}_q^d$. Consequently, using this relationship, we prove that the $L^2$ restriction conjectures for spheres hold in all dimensions when the test functions are restricted to homogeneous functions of degree zero.

\end{abstract}


\maketitle
\section{Introduction}
In recent years, several central open problems in Euclidean space have been actively studied in the finite field setting. Among these, the Kakeya problem, the restriction problem, and the Erd\H{o}s-Falconer distance problem have received significant attention and have been extensively investigated.\\

The Kakeya conjecture over finite fields, proposed by Wolff \cite{Wo99} in 1999, was completely resolved by Dvir \cite{Dv08} in 2008 using the polynomial method. His work introduced profound new ideas to the fields of harmonic analysis and discrete geometry (see, for example, \cite{Gu10, EOT10, GK15}).\\

The finite field restriction problem was initially posed and studied by Mockenhaupt and Tao \cite{MT04} in 2004 for various algebraic varieties over finite fields. 
Let $\mathbb{F}_q^d$, with $d \geq 2$, denote the $d$-dimensional vector space over a finite field $\mathbb{F}_q$ with $q$ elements. In the finite field setting, the paraboloid $P$ in $\mathbb{F}_q^d$ is defined as the set of points $x=(x_1,\ldots, x_d)\in \mathbb F_q^d$ satisfying the equation,  $x_d=x_1^2+\cdots+x_{d-1}^2:$
$$ P:=\{x\in \mathbb F_q^d: x_d=x_1^2+\cdots+x_{d-1}^2\}.$$ 
The work of Mockenhaupt and Tao on the finite field restriction problem primarily focused on the paraboloid and yielded good results in low dimensions. Their results on the paraboloid have been extended and enhanced to higher dimensions through the interest and efforts of many researchers (refer to \cite{IK09, LL10, Le13, SZ17, Le14, Ko20, IKL20, RS18, Le20} for example). On the other hand, the restriction problem for spheres over finite fields was first studied by Iosevich and the second listed author \cite{IK10}.
Compared to the paraboloid, the known results and techniques for spheres are relatively limited and have received less attention. This is unfortunate, as the restriction problem for spheres has proven highly useful when applied to other problems. For instance, as shown in \cite{CEHIK, KS10, CKY, CKY19}, restriction estimates for spheres have direct applications to the Erd\H{o}s-Falconer distance problem, which was initially posed and studied by Iosevich and Rudnev in \cite{IR06}.\\

The main goal of this paper is to introduce a new approach to the study of the restriction estimates for spheres and to engage readers' interest in this problem. As a main consequence of the approach, we show that the $L^2$ restriction conjectures for spheres are true in the specific case when the test functions are the homogeneous functions of degree zero. In addition, we address more reliable conjecture on the restriction problem for spheres over finite fields, which was not clearly stated in previous papers. \\

To precisely state our results in this introduction, we review the restriction problem for spheres in the finite field setting. In addition, we introduce a new operator $\mathcal{S}$, called the $\mathcal{S}$-operator, which is closely related to the restriction estimates for spheres over finite fields.\\

\subsection{The restriction operator $R_{S_j^{d-1}}$ for spheres $S_j^{d-1}$}
We begin with some notation and definitions. 
\begin{definition} We denote by $\mathcal{F}(\mathbb F_q^d \to \mathbb C)$ the set of all functions $g: \mathbb F_q^d \to \mathbb C.$ \end{definition}
We shall often use the notation $\mathcal{A}$ to indicate a subset of $\mathcal{F}(\mathbb F_q^d \to \mathbb C).$

Let $1\le p< \infty$ and let  $V, W$ be  sets.  Given functions $f: V\to \mathbb C$ and $g: W\to \mathbb C,$  we define
$$ ||f||_{L^p(V)}:= \left(\frac{1}{|V|} \sum_{x\in V} |f(x)|^p \right)^{\frac{1}{p}},    \quad ||g||_{L^\infty(V)}:=\max_{x\in V} |f(x)|,$$
and
$$||g||_{\ell^p(W)}:= \left( \sum_{m\in  W} |g(m)|^p \right)^{\frac{1}{p}},    \quad ||g||_{\ell^\infty(W)}:=\max_{m\in W} |g(m)|,$$
where $|V|$ denotes the cardinality of the set $V.$ Notice that $L^p(V)$ is associated with the normalizing counting measure on the set $V.$ On the other hand,  $\ell^p(W)$ is related to the counting measure supported on $W.$\\

For each finite field $\mathbb F_q$,  we choose a $j\in \mathbb F_q.$ 
Then we define a sphere $S_j^{d-1} \subset \mathbb F_q^d$  as the set
$$ S_j^{d-1}=\{x\in \mathbb F_q^d:  ||x||=j\}.$$
Here, and throughout the paper,   $||x||:=\sum\limits_{i=1}^d {x_i}^2$ for $x=(x_1, \ldots, x_d)\in \mathbb F_q^d.$
Let $g: \mathbb{F}_q^d \to \mathbb{C}$ be a function. Its Fourier transform, represented by $\widehat{g}$, is defined as  
\begin{equation}\label{defghat}
\widehat{g}(x) = \sum_{m \in \mathbb{F}_q^d} \chi(-m \cdot x) g(m),
\end{equation}
where $\chi$ is a fixed non-trivial additive character of $\mathbb F_q.$
The restriction operator $R_{S_j^{d-1}}$ for the sphere $S_j^{d-1}$ is defined by the relation
$$ R_{S_j^{d-1}}g(x):= \widehat{g}|_{S_j^{d-1}}(x).$$


Recall that  the additive character $\chi$  possesses the property of orthogonality: for any $a\in \mathbb F_q^*,$
$$ \sum_{t\in \mathbb F_q} \chi(at)=0.$$
Using the orthogonality of $\chi$, one can easily deduce 
the Fourier inversion theorem stated as follows: 
\begin{equation}\label{FI} g(m)=\frac{1}{q^d} \sum_{x\in \mathbb F_q^d} \chi(m\cdot x) \widehat{g}(x).\end{equation}
Using the orthogonality of $\chi$, it is not hard to deduce the Plancherel theorem:
$$ || \widehat{g} ||_{L^2(\mathbb F_q^d)} = ||g||_{\ell^2(\mathbb F_q^d)}.$$
In other words, the Plancherel theorem states that 
$$ \sum_{x\in \mathbb F_q^d} |\widehat{g}(x)|^2 = q^{d} \sum_{m\in \mathbb F_q^d} |g(m))|^2.$$

\begin{definition}[Restriction problem for spheres] Let $1\le p,r \le \infty.$ We define $R_{S_j^{d-1}}(p\to r)$ to be the smallest constant such that
the restriction estimate
\begin{equation}\label{RE}||\widehat{g}||_{L^r( S_j^{d-1})} \le R_{S_j^{d-1}}(p\to r) ||g||_{\ell^p(\mathbb F_q)}\end{equation} 
holds for all functions $g: \mathbb F_q^d \to \mathbb C.$  The restriction problem for the sphere $S_j^{d-1}$ is to determine all exponents $1\le p, r\le \infty$ such that $R_{S_j^{d-1}}(p\to r)\lesssim 1,$ where $A\lesssim B$ means that there exists a constant $C>0$ independent of $q$ such that  $A \le C B.$   We use $A\sim B$  to indicate that  $A\lesssim B$  and $B\lesssim A.$
We denote $R_{S_j^{d-1}}^\mathcal{A}(p\to r)\lesssim 1$ if the restriction estimate \eqref{RE} holds true for all functions $g\in \mathcal{A} \subseteq \mathcal{F}(\mathbb F_q^d \to \mathbb C).$
\end{definition}



\subsection{Spherical restriction conjectures and main results}
Throughout this paper, we will focus on finding all exponents $p$ such that $R_{S_j^{d-1}}(p\to 2)\lesssim 1$. This question is called the $L^2$ restriction problem for spheres over finite fields. Recall by the norm nesting property that if  
$R_{S_j^{d-1}}(p\to 2)\lesssim 1,$ then  $R_{S_j^{d-1}}(p_1\to 2) \lesssim 1$  for all $p_1$ with $1\le p_1 \le p \le \infty.$  Hence, the $L^2$ restriction problem for the sphere is to find the optimal $p$ value which is the largest number $p$ with $R_{S_j^{d-1}} (p\to 2) \lesssim 1.$ In two dimensions, the problem was completely solved by Iosevich and the second listed author in \cite{IKR}, showing that the optimal value of $p$ is $4/3.$  \\

For three and higher dimensions, $d\ge 3$,  they deduced by the Stein-Tomas argument that 
\begin{equation} \label{STR}
 R_{S_j^{d-1}}(p\to 2)\lesssim 1 \quad \mbox{if}~~ 1 \le p\le \frac{2d+2}{d+3}.
\end{equation}
In the case of Euclidean spaces, this range of $p$ cannot be improved, whereas in the case of finite fields, it is conjectured that it can be improved in special cases. 
However, precise statements of conjectures based on concrete examples are rarely found in existing papers.  For this reason,  we introduce a more refined conjecture and present partial evidence supporting the validity of our conjecture as the main result.\\

To address the conjecture, let us first recall that the quadratic character \( \eta \) of \( \mathbb{F}_q^* \) is the multiplicative character of \( \mathbb{F}_q^* \) that maps square numbers in \( \mathbb{F}_q^* \) to 1 and non-square numbers to -1:
\[
\eta(t) :=
\begin{cases}
1, & \text{if } t \text{ is a square number in } \mathbb{F}_q^*, \\
-1, & \text{if } t \text{ is a non-square number in } \mathbb{F}_q^*.
\end{cases}
\]
The value of $\eta(-1)$ depends on $\mathbb{F}_q$. More specifically, we have $\eta(-1) = -1$ if and only if $q \equiv 3 \pmod{4}$, and $\eta(-1) = 1$ if and only if $q \equiv 1 \pmod{4}$ (see Remark 5.13, \cite{LN97}).\\

We now introduce the $L^2$ restriction conjecture for spheres $S_j^{d-1}$. The background for this conjecture will be provided in the next section, where we deduce necessary conditions for $R_{S_j^{d-1}}(p\to r)\lesssim 1.$
\begin{conjecture} \label{conj1}  Suppose that $d\ge 3$  and $j\in \mathbb F_q^*.$ Then the following statements are true: 
\begin{enumerate}
\item If $d$ is even, then $R_{S_j^{d-1}}(p\to 2)\lesssim 1 $ if and only if $1\le p \le \frac{2d+4}{d+4}.$
\item If $d\equiv 1 \pmod{4}$ and $ \eta(j)=-1,$ then $R_{S_j^{d-1}}(p\to 2)\lesssim 1 $ if and only if $1\le p \le \frac{2d+6}{d+5}.$
\item If $d\equiv 1 \pmod{4}$ and $ \eta(j)=1,$ then $R_{S_j^{d-1}}(p\to 2)\lesssim 1 $ if and only if  $1\le p \le \frac{2d+2}{d+3}.$

\item If $d\equiv 3 \pmod{4}$ and $ \eta(-j)=-1,$ then $R_{S_j^{d-1}}(p\to 2)\lesssim 1 $ if and only if $1\le p \le \frac{2d+6}{d+5}.$

\item If  $d\equiv 3 \pmod{4}$ and $ \eta(-j)=1,$ 
then $R_{S_j^{d-1}}(p\to 2)\lesssim 1 $ if and only if  $1\le p \le \frac{2d+2}{d+3}.$

\end{enumerate}
\end{conjecture}

We have the following comments on Conjecture \ref{conj1}.
\begin{itemize}
\item In the conclusion of Conjecture \ref{conj1},  each necessary part for the boundedness of the $R_{S_j^{d-1}}(p\to 2)$  is true (see  Lemma \ref{lemNe} in Section \ref{sec2}). 

\item Hence, the known estimate \eqref{STR} implies that Conjecture \ref{conj1} holds for the  items $(3)$ and $(5)$.
\item In conclusion,  to settle the $L^2$ restriction conjectures for spheres over finite fields, it remains to establish the following.
\end{itemize}

\begin{conjecture}\label{conj2} Let $d\ge 3$ and $j\in \mathbb F_q^*.$ 
\begin{enumerate}
\item If $d$ is even, then $R_{S_j^{d-1}}(p\to 2)\lesssim 1 $ for  $1\le p \le \frac{2d+4}{d+4}.$
\item If $d\equiv 1 \pmod{4}$ and $ \eta(j)=-1$ (or  $d\equiv 3 \pmod{4}$ and $ \eta(-j)=-1$), then  we have
$R_{S_j^{d-1}}(p\to 2)\lesssim 1 $  for $1\le p \le \frac{2d+6}{d+5}.$
\end{enumerate}
    
\end{conjecture}



Proving or disproving Conjecture \ref{conj2} is considered a difficult problem. In fact, for any $d\ge 3,$ no result better than the Stein-Tomas result \eqref{STR} has been discovered, nor has any potential counterexample suggesting the conjecture might be false been identified. \\

We aim to provide some partial evidence that Conjecture \ref{conj2} is true.
As a key result for this purpose, we prove the weak version of the sharp  $L^2$ restriction estimates for spheres, where the test functions are restricted to the homogeneous functions on $\mathbb F_q^d$ of degree zero.  To precisely state our result,  let $\mathcal{H}:=\{g\in \mathcal{F}(\mathbb F_q^d\to \mathbb{C}):  g(m)=g(\lambda m) \quad \mbox{for all}~~ m\in \mathbb F_q^d,  \lambda\in \mathbb F_q^*\},$ which is a collection of the homogeneous functions of degree zero.
In \cite{KK14},  the authors showed that  the first part of Conjecture \ref{conj2} is true if the conclusion  $R_{S_j^{d-1}}(p\to 2) \lesssim 1$ is replaced  by a weaker condition that  $R^\mathcal{H}_{S_j^{d-1}}(p\to 2) \lesssim 1.$ As a major new result, stated below,  we will prove that the second statement of Conjecture \ref{conj2} also holds for such a weaker conclusion, $R^\mathcal{H}_{S_j^{d-1}}(p\to 2) \lesssim 1.$ 

\begin{theorem}\label{mainthm} If $d\equiv 1 \pmod{4}$ and $ \eta(j)=-1$ (or  $d\equiv 3 \pmod{4}$ and $ \eta(-j)=-1$), then  we have
$R^{\mathcal{H}}_{S_j^{d-1}}(p\to 2)\lesssim 1 $  for $1\le p \le \frac{2d+6}{d+5}.$
\end{theorem}

Let us highlight the important points of Theorem \ref{mainthm}.
\begin{itemize}
\item The range of $p$ in Conjecture \ref{conj2} (2) is much broader than that in Conjecture \ref{conj2} (1).  Hence, Proving our Theorem \ref{mainthm} for the special case of odd dimensions is distinct from proving the corresponding result for even dimensions and requires significantly more effort.

\item As we will see in Lemma \ref{Lem1.6}, a key ingredient in determining the values of \( p \) satisfying \( R_{S_j^{d-1}}^{\mathcal{H}}(p\to 2) \lesssim 1 \) is the \( L^2 \) restriction estimate for the \( j \)-homogeneous variety \( H_j^d \) in $\mathbb F_q^{d+1}$, which is defined in Definition \ref{Def1.5}.  

\item For example, after computing the sharp Fourier decay estimate on \( H_j^d \), 
one can apply the well-known Stein-Tomas argument to derive the \( L^2 \) restriction result for \( H_j^d \) 
and use it to find values of \( p \) satisfying \( R_{S_j^{d-1}}^{\mathcal{H}}(p\to 2) \lesssim 1 \). 
Using this method, it has been previously shown that \( R_{S_j^{d-1}}^{\mathcal{H}}(p\to 2) \lesssim 1 \) holds for all conjectured values of \( p \) satisfying 
\( R_{S_j^{d-1}}(p\to 2) \lesssim 1 \), except for cases $(2)$ and $(4)$ in Conjecture \ref{conj1}, which are the same as the conditions in our Theorem \ref{mainthm}.
However, under the conditions of our Theorem \ref{mainthm}, this approach leads to only very weaker result, $ R_{S_j^{d-1}}^{\mathcal{H}}(p\to 2) \lesssim 1$ for $1\le p \le \frac{2d+2}{d+3}.$ 

\item To prove Theorem \ref{mainthm}, we will also apply the result of the $L^2$ restriction estimate for $H_j^d$. 
However, to derive the $L^2$ restriction result, we will use a more efficient method than the well-known Stein-Tomas argument over finite fields. 
Roughly speaking, we derive a refined $L^2$ restriction estimate by explicitly computing the Fourier transform of the indicator function of $H_j^d$ 
when the test function is a characteristic function on $\mathbb F_q^{d+1}.$
Then, by applying the pigeonhole principle, we extend the result to general functions. 
We emphasize that our method can recover the sharp $L^2$ restriction result previously obtained using the Stein-Tomas method.
Indeed, in all cases, we will establish the sharp $L^2$ restriction result for $H_j^d$ using our approach (see Proposition \ref{Pro1.8}).

\end{itemize}

\subsection{Remaining part of this paper}
This paper will be organized as follows.
In Section \ref{sec2}, we introduce the necessary conditions for the boundedness of the restriction operator to the sphere $S_j^{d-1}$. Based on these conditions, we formulate the restriction conjecture for the sphere. In Section \ref{sec3}, we collect useful facts for the proof of Theorem \ref{mainthm} and provide its proof based on these facts.  In the final section, we complete the proof of the \( L^2 \) restriction results for the \( j \)-homogeneous variety \( H_j^d \subset \mathbb{F}_q^{d+1} \).

\section{Necessary conditions for the boundedness of  $R_{S_j^{d-1}}(p\to r)$}\label{sec2}
Some conjectures regarding restriction estimates for the spheres are mentioned in existing literature. However,  it seems that they have not been clearly stated. Therefore,  this section is devoted to presenting a more evidence-based conjecture.\\

It is well known that the main obstacle preventing the boundedness property that $R_{S_j^{d-1}}(p\to r)\lesssim 1$ arises when a large affine subspace $H$ is contained in the sphere $S_j^{d-1} \subset \mathbb F_q^d.$ Indeed,  assuming that $V$ contains an affine subspace $H$ of dimension $k$, namely $|H|=|\mathbb F|^k$, Mockenhaupt and Tao \cite{MT04} observed that  necessary conditions for the bound $R_{S_j^{d-1}}(p\to r) \lesssim 1$ can be given as follows:
\begin{equation}\label{nece1}  1\le p, r\le \infty, \quad
p'\ge \frac{2d}{d-1} \quad\mbox{and}\quad p'\ge \frac{r'(d-k)}{(r'-1)(d-1-k)}, 
\end{equation}
where $p'$ and $r'$ indicate H\"{o}lder's conjugates:
$$\frac{1}{p} + \frac{1}{p'}=1 
\quad \mbox{and} \quad  \frac{1}{r} + \frac{1}{r'}=1.$$
Notice that  the conditions \eqref{nece1} are the same as the following statement:
$(1/p, 1/r)$ lies on the convex hull of the four points
\begin{equation}\label{Boxnece} (1, 0),  (1,1),  \left(\frac{d+1}{2d},  1\right),   \left(\frac{d+1}{2d},  \frac{d^2-dk-d+k}{2d^2-2dk-2d}\right),\end{equation}
where $k$ denotes the dimension of an affine subspace $H$ lying on $S_j^{d-1},$ which means  $|H|=q^k.$
In particular, when $r=2,$  the necessary conditions \eqref{nece1} for $R_{S_j^{d-1}}(p\to 2) \lesssim 1$   become
$$  1\le p\le \infty, \quad
p'\ge \frac{2d}{d-1} \quad\mbox{and}\quad p'\ge \frac{2(d-k)}{(d-1-k)}. $$
In other words,  the necessary conditions for $R_{S_j^{d-1}}(p\to 2) \lesssim 1$  take the following:
\begin{equation}\label{r2conj} 1\le p \le \frac{2d^3-2d^2-2d^2k+2dk}{d^3-2d^2-d^2k +dk+d} .\end{equation}
The size of an affine subspace lying in $S_j^{d-1}$  is well known as follows (see Lemma 1.13 in \cite{KPV}):

\begin{lemma}\label{sizeASS} Let $S_j^{d-1}$ be the sphere in $\mathbb F_q^d$ with $j\ne 0.$ Then the following statements hold:
\begin{enumerate} \item If $d\ge 2$ is even , then $S_j^{d-1}$ contains an affine subspace $H$ with $|H|=q^{(d-2)/2}.$

\item If $d\equiv 1 \pmod{4}$ and $ \eta(j)=-1,$ then $S_j^{d-1}$ contains an affine subspace $H$ with $|H|=q^{(d-3)/2}.$
\item If $d\equiv 1 \pmod{4}$ and $ \eta(j)=1,$ then $S_j^{d-1}$ contains an affine subspace $H$ with $|H|=q^{(d-1)/2}.$

\item If If $d\equiv 3 \pmod{4}$ and $ \eta(-j)=-1,$  then  $S_j^{d-1}$ contains an affine subspace $H$ with $|H|=q^{(d-3)/2}.$
\item If If $d\equiv 3 \pmod{4}$ and $ \eta(-j)=1,$ then  $S_j^{d-1}$ contains an affine subspace $H$ with $|H|=q^{(d-1)/2}.$
\end{enumerate}
\end{lemma}

Combining \eqref{Boxnece} and Lemma \ref{sizeASS}, a direct computation yields necessary conditions on the boundedness of $R_{S_j^{d-1}}(p\to r),$ which can be conjectured as sufficient conditions. More precisely, one can obtain the following conjecture.

\begin{conjecture} Let $1\le p, r\le \infty.$ Then $R_{S_j^{d-1}}(p\to r) \lesssim 1$ provided that  each of the following statements holds:  
\begin{enumerate} \item  $d\ge 2$ is even and $(1/p, 1/r)$ lies in the convex hull of the points   $$(1, 0),  (1,1),  \left(\frac{d+1}{2d},  1\right),   \left(\frac{d+1}{2d},  \frac{d^2-d+2}{2d^2}\right),$$
\item $d\equiv 1 \pmod{4}, ~ \eta(j)=-1$ and $(1/p, 1/r)$ lies in the convex hull of the points  $$ (1, 0),  (1,1),  \left(\frac{d+1}{2d},  1\right),   \left(\frac{d+1}{2d},  \frac{d^2+3}{2d^2+2d}\right),$$
\item  $d\equiv 1 \pmod{4}, ~ \eta(j)=1$ and $(1/p, 1/r)$ lies in the convex hull of the points 
$$ (1, 0),  (1,1),  \left(\frac{d+1}{2d},  1\right),   \left(\frac{d+1}{2d},  \frac{d+1}{2d}\right),$$
\item  $d\equiv 3 \pmod{4}, ~ \eta(-j)=-1$  and $(1/p, 1/r)$ lies in the convex hull of the points   $$(1, 0),  (1,1),  \left(\frac{d+1}{2d},  1\right),   \left(\frac{d+1}{2d},  \frac{d^2+3}{2d^2+2d}\right),$$
\item  $d\equiv 3 \pmod{4},~ \eta(-j)=1$ and $(1/p, 1/r)$ lies in the convex hull of the points  
$$ (1, 0),  (1,1),  \left(\frac{d+1}{2d},  1\right),   \left(\frac{d+1}{2d},  \frac{d+1}{2d}\right).$$
\end{enumerate}
    
\end{conjecture}

In particular,  combining \eqref{r2conj} and Lemma \ref{sizeASS},  we obtain the following necessary conditions for $R_{S_j^{d-1}}(p\to 2) \lesssim 1$.  

\begin{lemma} \label{lemNe}
Suppose that $R_{S_j^{d-1}}(p\to 2) \lesssim 1$ for $1\le p \le \infty.$ Then the following statements hold true:
\begin{enumerate} \item If $d\ge 2$ is even, then  $1\le p \le \frac{2d+4}{d+4}.$

\item If $d\equiv 1 \pmod{4}$ and $\eta(j)=-1,$ then  $1\le p \le \frac{2d+6}{d+5}.$

\item If $d\equiv 1 \pmod{4}$ and $\eta(j)=1,$ then $1\le p \le \frac{2d+2}{d+3}.$

\item  If $d\equiv 3 \pmod{4}$ and $\eta(-j)=-1,$ then $1\le p \le \frac{2d+6}{d+5}.$

\item If  $d\equiv 3 \pmod{4}$  and $\eta(-j)=1,$ then $1\le p \le \frac{2d+2}{d+3}.$
\end{enumerate}
\end{lemma}

From this lemma, it is natural to consider Conjecture \ref{conj1}.

\section{Preliminary lemmas and the  proof of the main result, Theorem \ref{mainthm}} \label{sec3}
We first collect the key lemmas that play an important role in proving our main Theorem \ref{mainthm}. To this end, we begin by introducing a new operator which is called the $\mathcal{S}$ operator.

\begin{definition} [The $\mathcal{S}$-operator]
Given a function $g: \mathbb F_q^d \to C,$ we define the $\mathcal{S}$-operator by the relation
$$ \mathcal{S}g(m, m_{d+1}):= \frac{1}{q} \sum_{t\in \mathbb F_q^*} \chi(tm_{d+1}) g(tm),$$
where $m\in \mathbb F_q^d, m_{d+1}\in \mathbb F_q.$ 
\end{definition} 

Note that $\mathcal{S}g$ is a complex-valued function defined on $\mathbb F_q^{d+1}$ while $g$ is defined on $\mathbb F_q^d.$  

\begin{lemma} \label{lemF}
For a function $g:\mathbb F_q^d\to \mathbb C,$  the Fourier transform of the function $\mathcal{S}g: \mathbb F_q^{d+1}\to \mathbb C$, denoted by $\widehat{\mathcal{S}g},$ is given by
\begin{equation}\label{DG2} \widehat{\mathcal{S}g}(x, s)=\left\{\begin{array}{ll} \widehat{g}(x/ s)&\quad\mbox{if}~~s\ne 0\\
0&\quad\mbox{if}~~s= 0, \end{array} \right.\end{equation}
where  $x\in \mathbb F_q^{d}$ and  $s\in \mathbb F_q.$ 
\end{lemma}
\begin{proof}
Let $(x, s)\in \mathbb F_q^d\times \mathbb F_q=\mathbb F_q^{d+1}.$  By the definitions of the $\mathcal{S}$ and its Fourier transform, it follows that
\begin{align*} \widehat{\mathcal{S}g}(x,s)=& \sum_{(m, m_{d+1})\in \mathbb F_q^d\times \mathbb F_q} \chi(-x\cdot m -sm_{d+1}) \mathcal{S}g(m, m_{d+1})\\
=& \frac{1}{q} \sum_{m\in \mathbb F_q^d} \sum_{t\in \mathbb F_q^*}  \chi(-x\cdot m) g(tm)  \sum_{m_{d+1}\in \mathbb F_q} \chi((t-s)m_{d+1}). 
\end{align*}
If $s=0, $  the sum over $m_{d+1}\in \mathbb F_q$ vanishes by the orthogonality of $\chi$, since $t\in \mathbb F_q^*.$  Hence, $\widehat{\mathcal{S}g}(x,0) =0$ for all $x\in \mathbb F_q^d.$
On the other hand, if $s\in \mathbb F_q^*,$ then it follows by the orthogonality of $\chi$ that
$$\widehat{\mathcal{S}g}(x,s)=\frac{1}{q} \sum_{m\in \mathbb F_q^d} \chi(-x\cdot m) g(sm) q =\widehat{g}(x/s).$$
This completes the proof.
\end{proof}

We also need to consider a $j$-homogeneous variety in $\mathbb F_q^{d+1}$.
Given a finite field $\mathbb F_q,$ let $j$ be a fixed non-zero element in $\mathbb F_q.$ The number $j\in \mathbb F_q^*$ is always considered to be selected depending on $q.$
\begin{definition} \label{Def1.5} Let $j$ be a non-zero element in $\mathbb F_q.$ The $j$-homogeneous variety, denoted by $H_j^d$, is a variety lying in $\mathbb F_q^{d+1}$ and is defined by
$$ H_j^d=\{(x, x_{d+1})\in \mathbb F_q^d\times \mathbb F_q: ||x||=jx_{x+1}^2\}.$$
In addition, we define the dual variety of the $j$-homogeneous variety $H_j^d$, denoted by $(H_j^d)^*$, as
$$ (H_j^d)^*:= H_{j^{-1}}^d.$$
\end{definition}


We introduce one of the preliminary key lemmas to deduce our main theorem. 
\begin{lemma}\label{Lem1.6} Let $\mathcal{A}$ be a subset of $\mathcal{F}(\mathbb F_q^d \to \mathbb C)$ and let $1\le p, s, r\le \infty.$
Suppose that the following two estimates hold for all $g\in \mathcal{A}$:
\begin{itemize} \item[(i)]  $$||\widehat{\mathcal{S}g}||_{L^r(H_j^d)} ~\lesssim ~\|\mathcal{S}g\|_{\ell^s(\mathbb F_q^{d+1})}.$$ 

\item[(ii)] $$||\mathcal{S}g||_{\ell^s(\mathbb F_q^{d+1})} \lesssim ||g||_{\ell^p(\mathbb F_q^d)}.$$
Then we have
$$ R_{S_j^{d-1}}^{\mathcal{A}}(p\to r)\lesssim 1.$$ \end{itemize}
\end{lemma}

\begin{proof}
 By the assumptions $(i), (ii),$   to complete the proof, it suffices to show that  the following estimate
 \begin{equation}\label{eqk1}
|| \widehat{g}||_{L^r(S_j^{d-1})} \sim ||\widehat{\mathcal{S}g}||_{L^r(H_j^d)}
 \end{equation}
 holds for all functions $g \in  \mathcal{A}.$ 
Notice that for each $s\in \mathbb F_q^*,$ we have
 $$\widehat{g}(x)=\sum_{m\in \mathbb F_q^d} \chi(-sm\cdot x ) g(sm).$$ 
 It follows that
\begin{align*} \|\widehat{g}\|^r_{L^r(S_j^{d-1})}&=\frac{1}{|S_j^{d-1}|}\sum_{x\in  S_j^{d-1}} |\widehat{g}(x)|^r = \frac{1}{|S_j^{d-1}| (q-1)} \sum_{s\in \mathbb F_q^*} \sum_{x\in S_j^{d-1}} |\widehat{g}(x)|^r\\
&=\frac{1}{|S_j^{d-1}| (q-1)} \sum_{s\in \mathbb F_q^*} \sum_{x\in S_j^{d-1}}\left|\sum_{m\in \mathbb F_q^d} \chi(-sm\cdot x) g(sm)\right|^r.\end{align*}
Using a change of variables, $x\to x/s$ for a fixed $s\ne 0$ and  observing that  $|H_j^{d}|\sim q^{d}\sim q|S_j^{d-1}|,$ we get
\begin{equation*} \|\widehat{g}\|^r_{L^r(S_j^{d-1})}
\sim\frac{1}{|H_j^{d}|} \sum_{s\in \mathbb F_q^*} \sum_{x\in \mathbb F_q^d: ||x||=js^2} \left|\sum_{m\in \mathbb F_q^d} 
\chi(-m\cdot x) g(sm)\right|^r.\end{equation*}
For a fixed $s\ne 0,$ once again we use a change of variables,  $m \to m/s,$ and then
$$ \|\widehat{g}\|^r_{L^r(S_j^{d-1})}\sim \frac{1}{|H_j^{d}|} \sum_{s\in \mathbb F_q^*} \sum_{x\in \mathbb F_q^d: ||x||=js^2}|\widehat{g}(x/s)|^r.$$
By using \eqref{DG2} in Lemma \ref{lemF}, we obtain the desired consequence that 
$$\|\widehat{g}\|^r_{L^r(S_j^{d-1})}\sim\frac{1}{|H_j^{d}|} \sum_{(x, s)\in H_j^{d}} \left |\widehat{Sg}(x, s)\right|^r=\|\widehat{\mathcal{S}g}\|^r_{L^r(H_j^{d})}.$$
\end{proof}

Lemma \ref{Lem1.6} is useful in that it says that the restriction problems for the spheres $S_j^{d-1}$ in $d$-dimensions can be reduced to those for the homogeneous varieties $H_j^{d}$ in $(d+1)$-dimensions, respectively. 
Moreover,  from the literature, we see that the conjecture for the bound $R_{S_j^{d-1}}(p\to r)\lesssim 1$ is the same as that for the bound $R_{H_j^d}(p\to r)\lesssim 1$ (see \cite{Ko20, Le14, KLP22}). Similarly, one can also conjecture that the restriction estimates of the $S_j^{d-1}$ and $H_j^{d}$ are the same.
\subsection{Restriction estimates for the $j$-homogeneous varieties}
We need the $L^2$ restriction estimates for the $j$-homogeneous variety $H_j^d$ in $\mathbb F_q^{d+1}.$ These results will be combined with Lemma \ref{Lem1.6} in order to prove the main result, Theorem \ref{mainthm}. 

\begin{definition} For $1\le p, r\le \infty,$  we denote by $R_{H_j^{d}}(p\to r)$ the smallest constant such that 
the restriction estimate
$$ ||\widehat{G}||_{L^r(H_j^{d})} \le R_{H_j^d}(p\to r) ||G||_{\ell^p(\mathbb F_q^{d+1})}$$
holds true for all functions $G$ on $\mathbb F_q^{d+1}.$ The restriction problem for the $j$-homogeneous variety $H_j^{d}$ in $\mathbb F_q^{d+1}$ is to determine  all exponents $1\le p, r\le \infty$ such that  
$$R_{H_j^{d}}(p\to r) \lesssim 1.$$
\end{definition}

Now we states our results on the $L^2$ restriction estimates for the $j$-homogeneous variety $H_j^d$ in $\mathbb F_q^{d+1}.$
\begin{proposition}\label{Pro1.8}     For the $j$-homogeneous variety $H_j^d$, with $d\ge 2,$ in $\mathbb F_q^{d+1}$,  the following $L^2$ restriction estimates hold:
\begin{enumerate}
\item If $d$ is even, then $R_{H_j^{d}}(p\to 2)\lesssim 1 $ for all $1\le p \le \frac{2d+4}{d+4}.$

\item If $d\equiv 1 \pmod{4}$ and $ \eta(j)=-1,$ then $R_{H_j^{d}}(p\to 2)\lesssim 1 $ for all $1\le p \le \frac{2d+6}{d+5}.$

\item If $d\equiv 1 \pmod{4}$ and $\eta(j)=1,$ then $R_{H_j^{d}}(p\to 2)\lesssim 1 $ for all $1\le p \le \frac{2d+2}{d+3}.$

\item If $d\equiv 3 \pmod{4}$ and $ \eta(-j)=-1,$ then $R_{H_j^d}(p\to 2)\lesssim 1 $ for all $1\le p \le \frac{2d+6}{d+5}.$

\item If  $d\equiv 3 \pmod{4}$  and $\eta(-j)=1,$ then $R_{H_j^{d}}(p\to 2)\lesssim 1 $ for all $1\le p \le \frac{2d+2}{d+3}.$
\end{enumerate}
\end{proposition}

The proof of Proposition \ref{Pro1.8} will be presented in the next section (Section \ref{Resproof}).
For a moment, let us accept Proposition \ref{Pro1.8}, which will be invoked to establish the proof of our main theorem, Theorem \ref{mainthm}.

\subsection{Proof of Theorem \ref{mainthm}}
In this subsection, using Proposition \ref{Pro1.8} and Lemma \ref{Lem1.6},  we will prove that the values of \( p \) satisfying 
\[
R_{S_j^{d-1}}^{\mathcal{H}}(p\to 2) \lesssim 1
\]
coincide with the conjectured values of \( p \) required for 
\[
R(p\to 2) \lesssim 1
\]
as given in Conjecture \ref{conj1}. This fact has already been established in \cite{KK14}, except for the case of 
Theorem \ref{mainthm}. However, we will provide a unified proof for all other cases, including Theorem \ref{mainthm}, using our approach here.
In other words, we prove the following statement.

\begin{theorem}\label{thmj2} Let $d\ge 3$ and $j\in \mathbb F_q^*.$ 
\begin{itemize}
\item [(i)]If $d$ is even, then $R^{\mathcal{H}}_{S_j^{d-1}}(p\to 2)\lesssim 1 $ for  $1\le p \le \frac{2d+4}{d+4}.$
\item [(ii)] If $d\equiv 1 \pmod{4}$ and $ \eta(j)=-1$ (or  $d\equiv 3 \pmod{4}$ and $ \eta(-j)=-1$), then  we have
$R_{S_j^{d-1}}(p\to 2)\lesssim 1 $  for $1\le p \le \frac{2d+6}{d+5}.$

\item [(iii)] If $d\equiv 1 \pmod{4}$ and $ \eta(j)=1$ (or  $d\equiv 3 \pmod{4}$ and $ \eta(-j)=1$), then  we have
$R_{S_j^{d-1}}(p\to 2)\lesssim 1 $  for $1\le p \le \frac{2d+2}{d+3}.$
\end{itemize}
    
\end{theorem}

\begin{proof}
To complete the proof,  we first invoke Lemma \ref{Lem1.6} from which  taking $\mathcal{A}=\mathcal{H}, r=2,$ and $s=p,$  it is enough to show that
for all homogeneous functions $g$ of degree zero  in $\mathcal{H},$ we have
\begin{equation}\label{eqkm}
||\widehat{\mathcal{S}g}||_{L^2(H_j^d)} ~\lesssim ~\|\mathcal{S}g\|_{\ell^p(\mathbb F_q^{d+1})}.
\end{equation}

and
\begin{equation}\label{eqkmm}
||\mathcal{S}g||_{\ell^p(\mathbb F_q^{d+1})} \lesssim ||g||_{\ell^p(\mathbb F_q^d)},
\end{equation}
where the index `$p$' is the value corresponding to each statement in Theorem  \ref{thmj2}.
The inequality \eqref{eqkm} follows immediately from Proposition \ref{Pro1.8}. 
To prove the inequality \eqref{eqkmm},  let $g\in \mathcal{H}$. Then $g(tm)=g(m)$  for all nonzero $t\in \mathbb F_q^*$ and $m$ in $\mathbb F_q^d.$ Hence, for $(m, m_{d+1})\in \mathbb F_q^d\times \mathbb F_q,$  we see that

$$\mathcal{S}g(m, m_{d+1})= \frac{1}{q} \sum_{t\in \mathbb F_q^*} \chi(tm_{d+1}) g(m)= \frac{g(m)}{q} (\delta_0(m_{d+1})-1),$$
where $\delta_0(m_{d+1})=1$ for $m_{d+1}=0,$  and $0$ otherwise.
It therefore follows that
\begin{align*}    
||\mathcal{S}g||^p_{\ell^p(\mathbb F_q^{d+1})} &=\sum_{(m, m_{d+1})\in \mathbb F_q^d\times \mathbb F_q} \left| \frac{g(m)}{q} ( q\delta_0(m_{d+1}) -1) \right|^p\\
&=\sum_{m\in \mathbb F_q^d, m_{d+1}\in \mathbb F_q^*}\left| \frac{g(m)}{q} ( -1) \right|^p  + \sum_{m\in \mathbb F_q^d}  \left| \frac{g(m)}{q} ( q -1) \right|^p\\
&= \frac{q-1}{q^p} \|g\|_{\ell^p(\mathbb F_q^{d})}^p + \left(\frac{q-1}{q}\right)^p \|g\|_{\ell^p(\mathbb F_q^{d})}^p  \le 2 \|g\|_{\ell^p(\mathbb F_q^{d})}^p.\end{align*} 
This proves the inequality \eqref{eqkmm} and thus the proof of Theorem \ref{thmj2} is complete. \end{proof}

\section{Proof of Proposition \ref{Pro1.8} (Restriction results on $H_j^d$)} \label{Resproof}
In this section, we provide a  complete proof of  Proposition \ref{Pro1.8}, the result on the $L^2$ restriction estimates for the $j$-homogeneous variety $H_j^d.$ Before proceeding with the proof, let us introduce Gauss sums and some useful facts derived from them.

\subsection{ Gauss sums and their applications}
For $a\in \mathbb F_q^*$,  the Gauss sum $\mathcal{G}_a$  is defined by
$$ \mathcal{G}_a:=\sum_{s\in \mathbb F_q^*} \eta(s) \chi(as),$$
where $\eta$ denotes the quadratic character of $\mathbb F_q$.  \\

Now, we collect the basic facts related to the Gauss sum. 
Using a change of variables and a property of the quadratic character $\eta,$ we easily see that for $a\in \mathbb F_q^*,$
$$ \sum_{s\in \mathbb F_q} \chi(as^2)= \mathcal{G}_a= \eta(a) \mathcal{G}_1.$$
More generally, completing a square,  it follows that  for $a\in \mathbb F_q^*$ and $b\in \mathbb F_q,$ 
\begin{equation} \label{CS2}\sum_{s\in \mathbb F_q} \chi(as^2+bs) = \eta(a) \mathcal{G}_1 \chi\left(-\frac{b^2}{4a}\right).\end{equation}
The modulus of the Gauss sum $\mathcal{G}_a$ is exactly the same as $\sqrt{q}$ (see \cite{LN97}).  It is clear that $\mathcal{G}_a \overline{\mathcal{G}_a}=|\mathcal{G}_a|^2= q$ for any $a\in \mathbb F_q^*.$
Notice by a change of variables that 
$$ \overline{\mathcal{G}}_a = \sum_{s\in \mathbb F_q^*} \eta(s) \chi(-as) = \eta(-1)  \mathcal{G}_a.$$
Hence, we see that for any $a\in \mathbb F_q^*,$
\begin{equation}\label{ExplicitG2} \mathcal{G}_a^2=\eta(-1) q.\end{equation}
This equality is useful in computing   $\widehat{1_{H_j^d}},$ the Fourier transform of the indicator function of $j$-homogeneous variety $H_j^d.$ 
For simplicity,  we write
$$\widehat{H_j^d}:=\widehat{1_{H_j^d}},$$
where we identify the set $H_j^d$  with the indicator function of the set $H_j^d.$

\begin{lemma} \label{FF}  
For any \( M=(m, m_{d+1}) \in \mathbb{F}_q^d \times \mathbb{F}_q = \mathbb{F}_q^{d+1} \),  
the Fourier transform on the \( j \)-homogeneous variety \( H_j^d \) is given explicitly as follows:

\begin{enumerate}
\item If \( d \) is even:
\[
\widehat{H_j^d}(M) =
\begin{cases} 
q^d \delta_0(M) & \text{if } M \in (H_j^d)^*, \\
q^{\frac{d}{2}} (\eta(-1))^{(d+2)/2} \eta(j||m|| - m_{d+1}^2) & \text{if } M \notin (H_j^d)^*.
\end{cases}
\]

\item If \( d \equiv 1 \pmod{4} \) and \( \eta(j)=-1 \):
\[
\widehat{H_j^d}(M) =
\begin{cases} 
q^d \delta_0(M) - (1 - q^{-1}) q^{\frac{d+1}{2}} & \text{if } M \in (H_j^d)^*, \\
q^{\frac{d-1}{2}} & \text{if } M \notin (H_j^d)^*.
\end{cases}
\]

\item If \( d \equiv 1 \pmod{4} \) and \( \eta(j)=1 \):
\[
\widehat{H_j^d}(M) =
\begin{cases} 
q^d \delta_0(M) + (1 - q^{-1}) q^{\frac{d+1}{2}} & \text{if } M \in (H_j^d)^*, \\
- q^{\frac{d-1}{2}} & \text{if } M \notin (H_j^d)^*.
\end{cases}
\]

\item If \( d \equiv 3 \pmod{4} \) and \( \eta(-j)=-1 \):
\[
\widehat{H_j^d}(M) =
\begin{cases} 
q^d \delta_0(M) - (1 - q^{-1}) q^{\frac{d+1}{2}} & \text{if } M \in (H_j^d)^*, \\
q^{\frac{d-1}{2}} & \text{if } M \notin (H_j^d)^*.
\end{cases}
\]

\item If \( d \equiv 3 \pmod{4} \) and \( \eta(-j)=1 \):
\[
\widehat{H_j^d}(M) =
\begin{cases} 
q^d \delta_0(M) + (1 - q^{-1}) q^{\frac{d+1}{2}} & \text{if } M \in (H_j^d)^*, \\
- q^{\frac{d-1}{2}} & \text{if } M \notin (H_j^d)^*.
\end{cases}
\]
\end{enumerate}

\end{lemma}

\begin{proof}  We start by observing that
$$ \widehat{H_j^d}(M)= \sum_{X\in \mathbb F_q^{d+1}: x_1^2+ \cdots+x_d^2-jx_{d+1}^2=0} \chi(-M\cdot X),$$
Here, and throughout the paper, we write   $M=(m_1, \ldots, m_d, m_{d+1}),~ X=(x_1, \ldots, x_d, x_{d+1}) \in \mathbb F_q^{d+1}.$  By the orthogonality of $\chi$,  we can write
$$ \widehat{H_j^d}(M)= q^d \delta_{0}(M) + q^{-1} \sum_{t\in \mathbb F_q^*}  \left (\prod_{\ell=1}^d \sum_{x_{\ell}\in \mathbb F_q} \chi(t x_{\ell}^2-m_{\ell}x_{\ell})\right) \left(\sum_{x_{d+1}\in \mathbb F_q} \chi(-jtx_{d+1}^2-m_{d+1}x_{d+1})\right). $$
Invoking the inequality \eqref{CS2},  we obtain that
\begin{equation} \label{FT0} \widehat{H_j^d}(M) =q^d \delta_{0}(M) + q^{-1} \mathcal{G}^{d+1}_1 \eta(-j) \sum_{t\in \mathbb F_q^*}  \eta^{d+1}(t)  \chi\left( \frac{m_1^2+\cdots + m_d^2 -j^{-1} m_{d+1}^2}{-4t}\right).
\end{equation}

\textbf{(Case 1)} Let us prove the first part (1). Suppose that $d$ is even. Then $\eta^{d+1}=\eta.$  It follows from \eqref{FT0} that
\begin{equation} \label{FT1}\widehat{H_j^d}(M) =q^d \delta_{0}(M) + q^{-1} \mathcal{G}^{d+1}_1 \eta(-j) \sum_{t\in \mathbb F_q^*}  \eta(t)  \chi\left( \frac{m_1^2+\cdots + m_d^2 -j^{-1} m_{d+1}^2}{-4t}\right). \end{equation}
If $M\in (H_j^d)^*,$ then the statement (1)  is obvious since $\chi(0)=1$ and $\sum_{t\in \mathbb F_q^*} \eta(t)=0.$
Hence, we assume that $M\notin (H_j^d)^*.$ Then it is clear that $\delta_0(M)=0.$
By a change of variables and the properties of the quadratic character $\eta$,  one can easily observe that
$$ \sum_{t\in \mathbb F_q^*}  \eta(t)  \chi\left( \frac{m_1^2+\cdots + m_d^2 -j^{-1} m_{d+1}^2}{-4t}\right)=\eta(-m_1^2-\cdots-m_d^2+ j^{-1} m_{d+1}^2) \mathcal{G}_1.$$
In addition,  observe from \eqref{ExplicitG2} that
$$ \mathcal{G}_1^{d+2} = q^{\frac{d+2}{2}}  (\eta(-1))^{(d+2)/2}.$$
Then the statement (1) follows by combining those observations and  the inequality \eqref{FT1}. \\

\textbf{(Case 2)}  Let us prove the statements $(2)-(5).$  Suppose that $d$ is odd. Then $\eta^{d+1} \equiv 1$ and so  the inequality \eqref{FT0} becomes
\begin{align*} \widehat{H_j^d}(M) &=q^d \delta_{0}(M) + q^{-1} \mathcal{G}^{d+1}_1 \eta(-j) \sum_{t\in \mathbb F_q^*}  \chi\left( \frac{m_1^2+\cdots + m_d^2 -j^{-1} m_{d+1}^2}{-4t}\right)\\
&=q^d \delta_{0}(M) + q^{-1} \mathcal{G}^{d+1}_1 \eta(-j) \left( q \delta_0(m_1^2+ \cdots + m_d^2-j^{-1}m_{d+1}^2)-1\right),
\end{align*}
where the last equality follows by the orthogonality of $\chi.$
From \eqref{ExplicitG2},  it is not hard to see that
$$\mathcal{G}^{d+1}_1 \eta(-j)=q^{\frac{d+1}{2}}  (\eta(-1))^{(d+3)/2} \eta(j).$$
Therefore, we obtain that
$$ \widehat{H_j^d}(M) =q^d \delta_{0}(M) +  q^{\frac{d-1}{2}}  (\eta(-1))^{(d+3)/2} \eta(j) \left( q \delta_0(m_1^2+ \cdots + m_d^2-j^{-1}m_{d+1}^2)-1\right).$$
Finally, the statements (2)-(5) follow from the definition of $\delta_0$ and the following simple observations, respectively:
\begin{itemize} 
\item [(i)] If $d\equiv 1 \pmod{4}$ and $ \eta(j)=-1,$  then
$$(\eta(-1))^{(d+3)/2} \eta(j)=-1.$$
\item [(ii)] If $d\equiv 1 \pmod{4}$ and $ \eta(j)=1,$ then
$$(\eta(-1))^{(d+3)/2} \eta(j)=1.$$
\item [(iii)] If $d\equiv 3 \pmod{4}$ and $ \eta(-j)=-1,$ then
$$(\eta(-1))^{(d+3)/2} \eta(j)=-1.$$
\item [(iv)] If $d\equiv 3 \pmod{4}$ and $ \eta(-j)=1,$ then
$$(\eta(-1))^{(d+3)/2} \eta(j)=1.$$
\end{itemize}
\end{proof}

Since $ \widehat{H_j^d}(0,\ldots, 0) = |H_j^d|$ and $(0,\ldots, 0)\in (H_j^d)^*,$ it is a direct consequence of Proposition \ref{FF} that for $d\ge 2,$
$$ |H_j^d|\sim q^d.$$

Based on the explicit form of the Fourier transform on the $j$-homogeneous variety $H_j^d$  lying in $\mathbb F_q^{d+1},$  we are able to establish the following restricted strong type $L^2$ estimate for $H_j^d.$
\begin{proposition}\label{ProFF} Let $E\subset \mathbb F_q^{d+1}.$ Then the following statements hold:
\begin{enumerate}
\item If $d$ is even, then 
\[
\lVert \widehat{E} \rVert_{L^2(H_j^d)} \lesssim
\begin{cases} 
q^{\frac{1}{2}} \lvert E \rvert^{\frac{1}{2}} & \text{if } q^{\frac{d+2}{2}} \leq \lvert E\rvert \leq q^{d+1}, \\
q^{-\frac{d}{4}} \lvert E \rvert & \text{if } q^{\frac{d}{2}} \leq \lvert E \rvert \leq q^{\frac{d+2}{2}}, \\
\lvert E \rvert^{\frac{1}{2}} & \text{if } 1 \leq \lvert E \rvert \leq q^{\frac{d}{2}}.
\end{cases}
\]

\item If $d\equiv 1 \pmod{4}$ and $ \eta(j)=-1,$ then 
\[
\lVert \widehat{E} \rVert_{L^2(H_j^d)} \lesssim
\begin{cases} 
q^{\frac{1}{2}} \lvert E \rvert^{\frac{1}{2}} & \text{if } q^{\frac{d+3}{2}} \leq \lvert E\rvert \leq q^{d+1}, \\
q^{-\frac{d+1}{4}} \lvert E \rvert & \text{if } q^{\frac{d+1}{2}} \leq \lvert E \rvert \leq q^{\frac{d+3}{2}}, \\
\lvert E \rvert^{\frac{1}{2}} & \text{if } 1 \leq \lvert E \rvert \leq q^{\frac{d+1}{2}}.
\end{cases}
\]

\item If $d\equiv 1 \pmod{4}$ and $ \eta(j)=1,$ then 
\[
\lVert \widehat{E} \rVert_{L^2(H_j^d)} \lesssim
\begin{cases} 
q^{\frac{1}{2}} \lvert E \rvert^{\frac{1}{2}} & \text{if } q^{\frac{d+1}{2}} \leq \lvert E\rvert \leq q^{d+1}, \\
q^{-\frac{d-1}{4}} \lvert E \rvert & \text{if } q^{\frac{d-1}{2}} \leq \lvert E \rvert \leq q^{\frac{d+1}{2}}, \\
\lvert E \rvert^{\frac{1}{2}} & \text{if } 1 \leq \lvert E \rvert \leq q^{\frac{d-1}{2}}.
\end{cases}
\]

\item If $d\equiv 3 \pmod{4}$ and $ \eta(-j)=-1,$ then 
\[
\lVert \widehat{E} \rVert_{L^2(H_j^d)} \lesssim
\begin{cases} 
q^{\frac{1}{2}} \lvert E \rvert^{\frac{1}{2}} & \text{if } q^{\frac{d+3}{2}} \leq \lvert E\rvert \leq q^{d+1}, \\
q^{-\frac{d+1}{4}} \lvert E \rvert & \text{if } q^{\frac{d+1}{2}} \leq \lvert E \rvert \leq q^{\frac{d+3}{2}}, \\
\lvert E \rvert^{\frac{1}{2}} & \text{if } 1 \leq \lvert E \rvert \leq q^{\frac{d+1}{2}}.
\end{cases}
\]

\item If $d\equiv 3\pmod{4}$ and $ \eta(-j)=1,$ then 
\[
\lVert \widehat{E} \rVert_{L^2(H_j^d)} \lesssim
\begin{cases} 
q^{\frac{1}{2}} \lvert E \rvert^{\frac{1}{2}} & \text{if } q^{\frac{d+1}{2}} \leq \lvert E\rvert \leq q^{d+1}, \\
q^{-\frac{d-1}{4}} \lvert E \rvert & \text{if } q^{\frac{d-1}{2}} \leq \lvert E \rvert \leq q^{\frac{d+1}{2}}, \\
\lvert E \rvert^{\frac{1}{2}} & \text{if } 1 \leq \lvert E \rvert \leq q^{\frac{d-1}{2}}.
\end{cases}
\]
\end{enumerate}
\end{proposition}

\begin{proof}
By a direct computation, to complete the proof,  it will be enough to establish all the following statements:

\begin{itemize}
\item [(i)] If $d$ is even, then  
$$ \|\widehat{E}\|_{L^2(H_j^d)} \lesssim \min\left\{ |E|^{\frac{1}{2}} + q^{-\frac{d}{4}}|E|, ~~ q^{\frac{1}{2}} |E|^{\frac{1}{2}}\right\}. $$
\item [(ii)] If the assumption of Proposition \ref{ProFF} (2)  or Proposition \ref{ProFF} (4)  holds, then
$$ \|\widehat{E}\|_{L^2(H_j^d)} \lesssim \min\left\{ |E|^{\frac{1}{2}} + q^{-\frac{d+1}{4}}|E|, ~~ q^{\frac{1}{2}} |E|^{\frac{1}{2}}\right\}. $$
\item [(iii)] If the assumption of Proposition \ref{ProFF} (3)  or Proposition \ref{ProFF} (5)  holds, then
$$ \|\widehat{E}\|_{L^2(H_j^d)} \lesssim \min\left\{ |E|^{\frac{1}{2}} + q^{-\frac{d-1}{4}} |E|, ~~ q^{\frac{1}{2}} |E|^{\frac{1}{2}}\right\}. $$
\end{itemize}
Indeed, it is not hard to notice that statement (i) leads directly to Proposition \ref{ProFF} (1), statement (ii) corresponds to Proposition \ref{ProFF} (2) and (4), and statement (iii) supports Proposition \ref{ProFF} (3) and (5).

Since $|H_j^d|\sim q^d$ for $d\ge 2,$  we have
$$ \|\widehat{E}\|_{L^2(H_j^d)} \sim  \left( q^{-d} \sum_{M \in H_j^d} |\widehat{E}(M)|^2 \right)^{1/2}.$$
For a simple notation,  let us define that 
$$ \Omega(E):= \sum_{M \in H_j^d} |\widehat{E}(M)|^2.$$

Then, to establish statements (i), (ii), (iii),  it suffices to prove  the following claims (i)', (ii)', and (iii)', respectively:
\begin{itemize}
\item [(i)'] If $d$ is even, then  
$$ \Omega(E) \lesssim \min\left\{ q^d|E| + q^{\frac{d}{2}}|E|^2, ~~ q^{d+1} |E|\right\}. $$
\item [(ii)'] If the assumption of Proposition \ref{ProFF}(2)  or Proposition \ref{ProFF}(4)  holds, then
$$ \Omega(E) \lesssim \min\left\{ q^d|E| + q^{\frac{d-1}{2}}|E|^2, ~~ q^{d+1} |E|\right\}. $$
\item [(iii)'] If the assumption of Proposition \ref{ProFF}(3)  or Proposition \ref{ProFF}(5)  holds, then
$$ \Omega(E) \lesssim \min\left\{ q^d|E| + q^{\frac{d+1}{2}}|E|^2, ~~ q^{d+1} |E|\right\}. $$
\end{itemize}

Now let us estimate $\Omega(E).$ By Plancherel's theorem,  it is clear that
\begin{equation} \label{FCE}\Omega(E)\le \sum_{M\in \mathbb F_q^{d+1}} |\widehat{E}(M)|^2 = q^{d+1} |E|.\end{equation}
By the definition of the Fourier transform,  we see that
$$ \Omega(E) = \sum_{X, Y\in E} \widehat{H_j^d}(X-Y).$$
We compute the sum in $X, Y\in E$ by separating the cases where the $X-Y$ is included in $(H_j^d)^*$ and where it is not. Namely, we write
\begin{equation} \label{OSE}\Omega(E) = \sum_{X, Y\in E : X-Y\in (H_j^d)^*} \widehat{H_j^d}(X-Y)~ + \sum_{X, Y\in E : X-Y\notin (H_j^d)^*} \widehat{H_j^d}(X-Y).\end{equation}

\textbf{(Proof of (i)')}  Let $d$ be even.  By the estimate \eqref{FCE}, we only need to show that 
$$\Omega(E) \lesssim  q^d|E| + q^{\frac{d}{2}}|E|^2. $$
Combining \eqref{OSE} and Lemma \ref{FF}(1),  we get the required result:
\begin{align*} \Omega(E)&= q^d\sum_{X, Y\in E : X-Y\in (H_j^d)^*} \delta_0(X-Y)~ + q^{\frac{d}{2}} (\eta(-1))^{(d+2)/2} \sum_{X, Y\in E : X-Y\notin (H_j^d)^*} \eta(j||x-y||-(x_{d+1}-y_{d+1})^2)\\
&\le q^d |E| + q^{\frac{d}{2}} |E|^2. \end{align*}

\textbf{(Proof of (ii)')} From the estimate \eqref{FCE}, it is enough to show that
$$\Omega(E)\lesssim q^d|E| + q^{\frac{d-1}{2}}|E|^2.$$

Assume that the hypothesis on Proposition \ref{ProFF}(2) holds, namely $d\equiv 1 \pmod{4}$ and $ \eta(j)=-1.$  Then we can invoke Lemma \ref{FF}(2). 
Combining \eqref{OSE} and Lemma \ref{FF}(2), it follows that
\begin{align*} \Omega(E)&= q^d\sum_{X, Y\in E : X-Y\in (H_j^d)^*} \delta_0(X-Y)~ -(1-q^{-1}) q^{\frac{d+1}{2}} \sum_{X, Y\in E : X-Y\in (H_j^d)^*} 1\\
+&  \sum_{X, Y\in E : X-Y\notin (H_j^d)^*} q^{\frac{d-1}{2}}.
 \end{align*}
Since the second term above is non-positive,  it is clear that
\begin{align*}  \Omega(E)&\le  q^d\sum_{X, Y\in E : X-Y\in (H_j^d)^*} \delta_0(X-Y)
+  \sum_{X, Y\in E : X-Y\notin (H_j^d)^*} q^{\frac{d-1}{2}}\\
&\le q^d |E| + q^{\frac{d-1}{2}} |E|^2.\end{align*}
This completes the proof of (ii)' under the assumption on Proposition \ref{ProFF}(2).  When the assumption on Proposition \ref{ProFF}(4) holds,  the proof can be similarly demonstrated by using Lemma \ref{FF}(4) instead of Lemma \ref{FF}(2).

\textbf{(Proof of (iii)')}  By using the estimate \eqref{FCE}, our task is only to show that 
$$\Omega(E)\lesssim q^d|E| + q^{\frac{d+1}{2}}|E|^2.$$ 
If the hypothesis of Proposition \ref{ProFF}(3) is satisfied, Lemma \ref{FF}(3) can be applied to \eqref{OSE}. Similarly, if the hypothesis of Proposition \ref{ProFF}(5) is satisfied, Lemma \ref{FF}(5) can be applied to \eqref{OSE}.
Hence, we see that
\begin{align*} \Omega(E)&= q^d\sum_{X, Y\in E : X-Y\in (H_j^d)^*} \delta_0(X-Y)~ +(1-q^{-1}) q^{\frac{d+1}{2}} \sum_{X, Y\in E : X-Y\in (H_j^d)^*} 1\\
-&  \sum_{X, Y\in E : X-Y\notin (H_j^d)^*} q^{\frac{d-1}{2}}.
 \end{align*}
Since the third term above is  non-positive,  it is clear that
\begin{align*}  \Omega(E)&\le q^d\sum_{X, Y\in E : X-Y\in (H_j^d)^*} \delta_0(X-Y)~ +(1-q^{-1}) q^{\frac{d+1}{2}} \sum_{X, Y\in E : X-Y\in (H_j^d)^*} 1\\
&\lesssim  q^d |E| + q^{\frac{d+1}{2}} |E|^2.\end{align*}
This completes the proof of (iii)'.  We have finished Proposition \ref{ProFF}.
\end{proof}

\subsection{Reduction lemma}
In this subsection, we introduce  notation and a useful lemma necessary for the reduction in the proof of Proposition \ref{Pro1.8}. 

\begin{definition}
Let $F: \mathbb F_q^{d+1} \to [0, 1].$ 
\begin{itemize}
\item [(i)]
For each non-negative integer $i\ge 0,$  we denote
$$ F_i:= \{X\in \mathbb F_q^{d+1}: 2^{-i-1} < F(X) \le 2^{-i}\}.$$
\item [(ii)] The function $\widetilde{F}$ is a step function on $\mathbb F_q^{d+1}$  defined as
$$ \widetilde{F}(X):= \sum_{i=0}^\infty  2^{-i} 1_{F_i}(X).$$
\end{itemize}
\end{definition}

Notice that the collection \(\{F_i\}\) consists of mutually disjoint sets.
It is clear that $\bigcup\limits_{i=0}^\infty F_i \subseteq \mathbb F_q^{d+1},$ which  implies that  $|F_i|\le q^{d+1}$ for all $i\ge 0.$
A simple but important observation is as follows. For each \( X \in \mathbb{F}_q^{d+1} \), we have  
\begin{equation} \label{FFt}
F(X) \leq \widetilde{F}(X) \leq 2F(X).\end{equation}

The following facts are necessary for the proof of Proposition \ref{Pro1.8}.
\begin{lemma}\label{RLem} Let $F: \mathbb F_q^{d+1} \to [0, 1]$ and let $1\le p < \infty.$ 
\begin{enumerate}
    \item  If $\sum\limits_{X\in \mathbb F_q^{d+1}} (F(x))^p =1,$  then $\sum\limits_{k=0}^\infty 2^{-pk} |F_k| \le 2^p.$

    \item We have $ ||\widehat{F}||_{L^2(H_j^d)} \sim ||\widehat{\widetilde{F}}||_{L^2(H_j^d)}.$
\end{enumerate}
\end{lemma}

\begin{proof} First,  let us prove Lemma \ref{RLem} (1).   By combining  \eqref{FFt} and  the assumption that $\sum\limits_{X\in \mathbb F_q^{d+1}} (F(x))^p =1,$   we see that
$$ \sum_{X\in \mathbb F_q^{d+1}} (\widetilde{F}(X))^p  \le 2^p.$$
It follows by the definition of $\widetilde{F}$ that
$$ \sum_{X\in \mathbb F_q^{d+1}} (\widetilde{F}(X))^p =\sum_{k=0}^\infty  \sum_{X\in F_k}  \left( \sum_{i=0}^\infty  2^{-i} 1_{F_i}(X)\right)^p = \sum_{k=0}^\infty 2^{-pk} |F_k|.$$ 
This proves the first part of Lemma \ref{RLem}.\\

Now we prove Lemma \ref{RLem} (2). It is enough to prove that
\[ \sum_{X \in H_j^{d}} |\widehat{F}(X)|^2 \sim
\sum_{X \in H_j^{d}} |\widehat{\widetilde{F}}(X)|^2 .
\]
By the definition of the Fourier transform, we have
\[
\sum_{X \in H_j^{d}} |\widehat{F}(X)|^2 = \sum_{X \in H_j^{d}} \left| \sum_{M \in \mathbb{F}^{d+1}} {F}(M) \chi(-X \cdot M) \right|^2
= \sum_{X \in H_j^{d}} \sum_{M, M' \in \mathbb{F}^{d+1}} {F}(M) {F}(M') \chi(-X \cdot (M - M')),
\]
where we use the fact that \({F}\) is real-valued. By the definition of the Fourier transform on \( H_j^{d} \), the above expression is equal to
\begin{equation}\label{LL}
\sum_{X \in H_j^{d}} |\widehat{{F}}(X)|^2 =  \sum_{M, M' \in \mathbb{F}^{d+1}} {F}(M) {F}(M') \widehat{H_j^d}(M - M').
\end{equation}
Similarly, we obtain
\begin{equation} \label{RR}
\sum_{X \in H_j^{d}} |\widehat{\widetilde{F}}(X)|^2 =  \sum_{M, M' \in \mathbb{F}^{d+1}} \widetilde{F}(M) \widetilde{F}(M') \widehat{H_j^d}(M - M').
\end{equation}
Since \(F \sim \widetilde{F}\) by \eqref{FFt}, it follows that
\[
\widetilde{F}(M) \widetilde{F}(M') \sim F(M) F(M') \quad \text{for all } M, M' \in \mathbb{F}^{d+1}.
\]
Moreover,  it follows from Lemma \ref{FF} that  the Fourier transform \( \widehat{H_j^d} \) can take at most three distinct real values. 
Thus, the right-hand side of equation \eqref{LL} is similar to that of equation \eqref{RR}, concluding the proof.
\end{proof}

\subsection{Complete Proof of Proposition \ref{Pro1.8}}

We will prove Proposition \ref{Pro1.8} by invoking the pigeonhole principle used in \cite{IKL20, KLP22, RS18}  together with our Proposition \ref{ProFF}.  For the convenience of readers, we state Proposition \ref{Pro1.8} in a more concise form and provide its proof below.

\begin{proposition} \label{FPK}  For the $j$-homogeneous variety $H_j^d$ in $\mathbb F_q^{d+1}$, where $d\ge 2,$    the following $L^2$ restriction estimates hold:
\begin{itemize}
\item [(i)]  If $d$ is even, then $R_{H_j^{d}}(p\to 2)\lesssim 1 $ for all $1\le p \le \frac{2d+4}{d+4}.$

\item [(ii)] If $d\equiv 1 \pmod{4}$ and $ \eta(j)=-1$ (or $d\equiv 3 \pmod{4}$ and $ \eta(-j)=-1$), then we have $R_{H_j^{d}}(p\to 2)\lesssim 1 $ for all $1\le p \le \frac{2d+6}{d+5}.$

\item [(iii)] If $d\equiv 1 \pmod{4}$ and $\eta(j)=1$ (or $d\equiv 3 \pmod{4}$  and $\eta(-j)=1$), then $R_{H_j^{d}}(p\to 2)\lesssim 1 $ for all $1\le p \le \frac{2d+2}{d+3}.$
\end{itemize}
\end{proposition}
\begin{proof} Before starting the proof, note that Proposition \ref{FPK} (i) matches Proposition \ref{Pro1.8} (1), Proposition \ref{FPK} (ii) implies Proposition \ref{Pro1.8} (2) and (4), and Proposition \ref{FPK} (iii) corresponds to Proposition \ref{Pro1.8} (3) and (5).

We will prove the conclusions (i), (ii), and (iii) of Proposition \ref{FPK} simultaneously.
To do this,  we first define $p=p(\alpha)$ for $ 0<\alpha < d$ as 
\begin{equation}\label{defp}
p = p(\alpha) := \frac{2\alpha+2}{\alpha+2}.
\end{equation}
For each statement of conclusions $(i), (ii)$ and $(iii)$, we correspondingly assign $\alpha=\frac{d}{2},  \frac{d+1}{2}, \frac{d-1}{2}$ so that the values of
$p=p(\alpha)=\frac{2d+4}{d+4}, \frac{2d+6}{d+5}, \frac{2d+2}{d+3}$  correspond to the critical value of $p$ for proving the statements $(i), (ii),$ and $(iii)$ of Proposition \ref{FPK}, respectively.  Therefore, in order to prove Proposition \ref{FPK},  for each  $p=p(\alpha), $   it suffices to show that 
the following restriction estimate holds for all functions $F$ on $\mathbb F_q^{d+1}:$
$$\| \widehat{F} \|_{L^2(H_j^d)} \lesssim \| F \|_{L^{p}\left(\mathbb F_q^{d+1}\right)} = 
\left( \sum_{X \in \mathbb F_q^{d+1} } |F(X)|^{p} \right)^{\frac{1}{p}}.$$
We now make reductions on the test functions $F$ on $\mathbb F_q^{d+1}.$
As usual,  we may assume that  the test function $F$ is a non-negative real-valued function on $\mathbb F_q^{d+1}.$
Furthermore, a normalization of the function $F$ allows us to assume that the test function $F$   satisfies the condition that 
\begin{equation*}\sum_{X \in \mathbb F_q^{d+1} } |F(X)|^{p} =1.\end{equation*}
Thus, our problem is reducing to showing that 
the inequality 
$$\| \widehat{F} \|_{L^2(H_j^d)} \lesssim 1 $$
holds for all real-valued functions $F: \mathbb F_q^{d+1} \to [0, 1]$ with the property that
\begin{equation*} \label{CK1}\sum_{X \in \mathbb F_q^{d+1} } (F(X))^{p} =1. \end{equation*}

By Lemma \ref{RLem} (1),  we also assume that 
\begin{equation}\label{Onp} \sum\limits_{k=0}^\infty 2^{-pk} |F_k| \lesssim 1,\end{equation}
where we recall that $F_k:= \{X \in \mathbb{F}_q^{d+1} : 2^{-k-1} < F(X) \le 2^{-k}\}.$
It is obvious from \eqref{Onp} that for all $k\ge 0,$
\begin{equation}\label{Onp1} |F_k| \lesssim 2^{pk}.\end{equation}

In summary,  our problem is reduced to showing that 
\begin{equation}\label{Goal1}\| \widehat{F} \|_{L^2(H_j^d)} \lesssim 1 \end{equation}
where the function $F$ satisfies the conditions \eqref{Onp} and \eqref{Onp1}. \\

By Lemma \ref{RLem} (2) and the definition of $\widetilde{F},$  we see that for any integer $N\ge 1,$
$$||\widehat{F}||_{L^2(H_j^d)} \sim ||\widehat{\widetilde{F}}||_{L^2(H_j^d)}  \le \sum_{k=0}^N 2^{-k} \|\widehat{{F_k}}\|_{L^2(H_j^d)}+ \sum_{k=N+1}^\infty 2^{-k} \|\widehat{{F_k}}\|_{L^2(H_j^d)}=: T(N) + R(N),$$
where we also used Minkowski's inequality. 
Our task is  to show that  $T(N)\lesssim 1$  and $R(N)\lesssim 1$  for some $N.$
We choose an integer $N$ such that 
$$N\ge (d+1) \log_2{q}.$$
Then,  since $\|\widehat{{F_k}}\|_{L^2(H_j^d)} \le |F_k| \le q^{d+1}$ for all $k\ge 0, $   we have
$$ R(N) \le q^{d+1} \sum_{k=N+1}^\infty 2^{-k} \le q^{d+1} 2^{-N} \le 1.$$

Finally, it remains to prove that 
$$ T(N)=\sum_{k=0}^N 2^{-k} \|\widehat{{F_k}}\|_{L^2(H_j^d)} \lesssim 1.$$

Let $[N]:=\{0, 1, \ldots, N\}.$ 
We partition $[N]$ into the following three parts:
\begin{align*}
 [N]_1 & := \{k\in [N]: 1\le 2^{pk} \le q^\alpha\},\\
 [N]_2 & := \{k\in [N]: q^\alpha \le 2^{pk} \le q^{\alpha+1}\},\text{ and}\\
 [N]_3 & := \{k\in [N]: q^{\alpha+1}\le 2^{pk} \le q^{d+1}\},
\end{align*}
where we recall that  $\alpha=\frac{d}{2},  \frac{d+1}{2}, \frac{d-1}{2}$ are the values defined for proving the statements $(i), (ii)$ and $(iii)$ of Proposition \ref{FPK}, respectively.
Then we can write 
\begin{align*}
T(N) & \le \sum_{k\in [N]_1} 2^{-k} \|\widehat{1_{F_k}}\|_{L^2(H_j^d)}+ 
\sum_{k\in [N]_2}2^{-k} \|\widehat{1_{F_k}}\|_{L^2(H_j^d)}+ 
\sum_{k\in [N]_3} 2^{-k} \|\widehat{1_{F_k}}\|_{L^2(H_j^d)}\\
&=: I + II +III
\end{align*}
To complete the proof,  our task is to  prove that   $I, II, III \lesssim 1.$

Each of the conclusions in Proposition \ref{ProFF} can be written in the following form: for $E\subset \mathbb F_q^{d+1},$
\begin{align} \label{Alp}
\lVert \widehat{E} \rVert_{L^2(H_j^d)} \lesssim
\begin{cases} 
q^{\frac{1}{2}} \lvert E \rvert^{\frac{1}{2}} & \text{if } q^{\alpha+1} \leq \lvert E\rvert \leq q^{d+1}, \\
q^{-\frac{\alpha}{2}} \lvert E \rvert & \text{if } q^{\alpha} \leq \lvert E \rvert \leq q^{\alpha+1}, \\
\lvert E \rvert^{\frac{1}{2}} & \text{if } 1 \leq \lvert E \rvert \leq q^{\alpha},
\end{cases}
\end{align}
where $\alpha= \frac{d}{2},  \frac{d+1}{2}, \frac{d-1}{2}$    correspond to Proposition \ref{ProFF} (1), Proposition \ref{ProFF} (2) and (4), and Proposition \ref{ProFF} (3) and (5), respectively.

Now we show that $I\lesssim 1.$
Since $ |F_k|\le 2^{pk}$ for all $k\in [N],$  it follows that
$$ I\le \sum_{k\in [N]_1 : 1\le |F_k|\le q^\alpha} 2^{-k} \|\widehat{1_{F_k}}\|_{L^2(H_j^d)} \lesssim \sum_{k\in [N]_1 : 1\le |F_k|\le q^\alpha} 2^{-k} |F_k|^{\frac{1}{2}} \le \sum_{k=0}^\infty  2^{-k} 2^{\frac{pk}{2}} \lesssim 1,$$
where the second inequality above follows from the estimate \eqref{Alp} and  the last inequality above is clear by the definition of $p$ given in \eqref{defp}. 

Next, let us prove $II \lesssim 1.$ It follows that
$$ II \le  \sum_{k\in [N]_2: 1\le |F_k|\le q^\alpha}2^{-k} \|\widehat{1_{F_k}}\|_{L^2(H_j^d)}
+ \sum_{k\in [N]_2: q^\alpha \le |F_k|\le q^{\alpha+1}} 2^{-k} \|\widehat{1_{F_k}}\|_{L^2(H_j^d)}$$
$$\lesssim \sum_{k\in [N]_2: 1\le |F_k|\le q^\alpha}2^{-k} |F_k|^{\frac{1}{2}}  +  \sum_{k\in [N]_2: q^\alpha \le |F_k|\le q^{\alpha+1}} 2^{-k} q^{-\frac{\alpha}{2}} |F_k| $$
$$\lesssim \sum_{k=0}^\infty  2^{-k} 2^{\frac{pk}{2}} + q^{-\frac{\alpha}{2}}   \sum_{k\in [N]_2}  2^{-k}  2^{pk}\lesssim   1+ q^{-\frac{\alpha}{2}} ~ 2^{(p-1)(\frac{\alpha+1}{p} \log_2{q}) }= 2,$$
where the last inequality follows by the observations that $p/2-1<0, p-1 >0,$ and 
    $ k\le  \frac{\alpha+1}{p} \log_2{q}$ for $k\in [N]_2.$
Finally, we show that $III\lesssim 1.$ By decomposing  the set $[N]_3$ into three parts,  we write the term $III$ as follows:
$$III\le \sum_{\substack{k\in [N]_3\\: 1\le |F_k|\le q^\alpha}} 2^{-k} \|\widehat{1_{F_k}}\|_{L^2(H_j^d)} + \sum_{\substack{k\in [N]_3\\: q^\alpha\le |F_k|\le q^{\alpha+1}}} 2^{-k} \|\widehat{1_{F_k}}\|_{L^2(H_j^d)} + \sum_{\substack{k\in [N]_3\\: q^{\alpha+1}\le |F_k|\le q^{d+1}}} 2^{-k} \|\widehat{1_{F_k}}\|_{L^2(H_j^d)}.$$

Applying the estimate \eqref{Alp} to each of the three items above and  using the property that  $|F_k|\le 2^{pk}$ for all $k\in [N],$  we obtain that
$$III\lesssim  \sum_{\substack{k\in [N]_3\\: 1\le |F_k|\le q^\alpha}} 2^{-k} |F_k|^{\frac{1}{2}} + \sum_{\substack{k\in [N]_3\\: q^\alpha\le |F_k|\le q^{\alpha+1}}} 2^{-k} q^{-\frac{\alpha}{2}} |F_k| + \sum_{\substack{k\in [N]_3\\: q^{\alpha+1}\le |F_k|\le q^{d+1}}} 2^{-k} q^{\frac{1}{2}} \lvert F_k \rvert^{\frac{1}{2}}$$
$$ \lesssim \sum_{k=0}^\infty  2^{-k} 2^{\frac{pk}{2}} +   \sum_{k\in [N]_3}  2^{-k} q^{-\frac{\alpha}{2}} q^{\alpha+1} +    q^{\frac{1}{2}}   \sum_{k\in [N]_3}   2^{-k} 2^{\frac{pk}{2}} =: III_1 + III_2+ III_3. $$
Now note from \eqref{defp} that $p/2-1 <0.$ Then it is obvious that $III_1 \sim 1.$
It remains to show that $III_2 + III_3 \lesssim 1.$ To do this,  we  observe that if $k\in [N]_3,$ then 
$$  \frac{\alpha+1}{p} \log_2{q} \le k. $$ 
In addition,  recall that  $ \frac{p}{2}-1 <0$. Then it follows that
$$ III_2 + III_3 \lesssim q^{\frac{\alpha}{2}+1}   \sum_{\frac{\alpha+1}{p} \log_2{q}\le k\le \infty}  2^{-k}  +    q^{\frac{1}{2}}   \sum_{\frac{\alpha+1}{p} \log_2{q}\le k\le \infty}   2^{-k} 2^{\frac{pk}{2}}$$
$$\lesssim q^{\frac{\alpha}{2}+1}~ 2^{- (\frac{\alpha+1}{p} \log_2{q})} +  q^{\frac{1}{2}} ~2^{(p/2-1)  (\frac{\alpha+1}{p} \log_2{q})} =1+1=2.$$
This completes the proof of Proposition \ref{Pro1.8}.  \end{proof}

\bibliographystyle{amsplain}

\begin{thebibliography}{10}


\bibitem{CEHIK} J. Chapman, M. B. Erdo\u{g}an, D. Hart,  A. Iosevich, and D. Koh, \emph{Pinned distance sets, k-simplices, Wolff's exponent in finite fields and sum-product estimates}, Math. Z. \textbf{271} (2012), no.1, 63-93.

\bibitem{CKY} D. Covert, Y. Pi, and D. Koh, \emph{On the sums of any $k$ points in finite fields,} SIAM J. Discrete Math, \textbf{30} (2016), no.1, 367-382. 

\bibitem{CKY19}D. Covert, Y. Pi, and D. Koh, \emph{The generalized $k$-resultant modulus set problem in finite fields,} J. Fourier Anal. Appl. \textbf{25} (2019), no.3, 1026-1052.

\bibitem{Dv08} Z. Dvir, \emph{On the size of Kakeya sets in finite fields}, J. Amer. Math. Soc. \textbf{22} (2009), no.4, 1093-1097.

\bibitem{EOT10} J. S. Ellenberg, R. Oberlin, and T. Tao, \emph{The Kakeya set and maximal conjectures for algebraic varieties over finite fields}, Mathematika, \textbf{56} (2010), no.1,  1-25.


\bibitem{Gu10} L. Guth, \emph{The endpoint case of the Bennett-Carbery-Tao multilinear Kakeya conjecture}, Acta Math. \textbf{205} (2010), no.2, 263-286.

\bibitem{GK15} L. Guth, N. H. Katz,\emph{On the Erd\H{o}s distinct distances problem in the plane}, Ann. of Math. (2) \textbf{181} (2015), no.1, 155-190.

\bibitem{IKR} A. Iosevich and D. Koh, \emph{Extension theorems for the Fourier transform associated with nondegenerate quadratic surfaces in vector spaces over finite fields}, Illinois J. of Mathematics, {\bf 52} (2008), no.2,  611-628.


\bibitem{IK09} A. Iosevich and D. Koh, \emph{Extension theorems for paraboloids in the finite field setting}, Math. Z. {\bf 266} (2010), no.2,  471-487.

\bibitem{IK10}
 A. Iosevich and D. Koh, \textit{Extension theorems for spheres in the finite field setting}, Forum. Math. \textbf{22} (2010), no.3, 457-483.

 
\bibitem{IKL20}
A. Iosevich, D. Koh and M. Lewko, \textit{Finite field restriction estimates for the paraboloid in high even dimensions},  J. Funct. Anal.  \textbf{278} (2020), no.11, 108450.
  
 
\bibitem{IR06} A. Iosevich and M. Rudnev, {\em Erd\H{o}s-Falconer distance problem in vector spaces over finite fields}, Trans. Amer. Math. Soc. \textbf{359} (2007), no. 12, 6127-6142.

\bibitem{KK14} H. Kang and D. Koh, \emph{Weak version of restriction estimates for spheres and paraboloids in finite fields}, J. Math. Anal. Appl. \textbf{419} (2014), no.2, 783-795. 

\bibitem{Ko20} D. Koh, \textit{Conjecture and improved extension theorems for paraboloids in the finite field setting}, Math. Z. \textbf{294} (2020),  no. 1-2, 51-69. 

\bibitem{KLP22} D. Koh, S. Lee, and T. Pham, \emph{On the cone restriction conjecture in four dimensions and applications in incidence geometry},  Int. Math. Res. Not.  (2022), no. 21, 17079-17111. 

\bibitem{KPV} D. Koh, T. Pham and L.A. Vinh, \emph{Extension theorems and a connection to the Erd\H{o}s-Falconer distance problem over finite fields}, J. Funct. Anal. \textbf{281} (2021), no.8, Paper No. 109137, 54 pp.

\bibitem{KS10} D. Koh and C. Shen, \emph {Sharp extension theorems and Falconer distance problems for algebraic curves in two dimensional vector spaces over finite fields}, Revista Matematica Iberoamericana, \textbf{28} (2012), no.1, 157-178.


\bibitem{Le13}  M. Lewko, \emph{New restriction estimates for the 3-d paraboloid over finite fields}, Adv. Math. {\bf 270} (2015), no.1, 457-479.

\bibitem{Le14}  M. Lewko, \emph{Finite field restriction estimates based on Kakeya maximal operator estimates},  J. Eur. Math. Soc. \textbf{21}  (2019),  no.12, 3649-3707. 

\bibitem{Le20} M. Lewko, \emph{Counting rectangles and an  improved restriction estimate for the paraboloid in $\mathbb F_p^3$,}  Proc. Amer. Math. Soc. \textbf{148} (2020), no.4,  1535-1543.

\bibitem{LL10} A.  Lewko and M. Lewko, \emph{Endpoint restriction estimates for the paraboloid over finite fields}, Proc. Amer. Math. Soc. \textbf{140}  (2012), no.6,  2013-2028.

\bibitem{LN97} R. Lidl and H. Niederreiter, \emph{Finite fields,} Cambridge University Press, (1997).

\bibitem{MT04} G. Mockenhaupt, and T. Tao, \emph{Restriction and Kakeya phenomena for finite fields}, Duke Math. J. {\bf 121} (2004), no.1, 35-74.

\bibitem{RS18} M. Rudnev and I. D. Shkredov, \emph{On the restriction problem for discrete paraboloid in lower dimension},  Adv. Math. \textbf{339} (2018), 657-671.

\bibitem{SZ17} S. Stevens and F. de Zeeuw, \emph{An improved point-line incidence bound over arbitrary fields,}  Bull. Lond. Math. Soc. \textbf{49} (2017),  no.5, 842-858.

\bibitem{Wo99} T. Wolff, \emph{Recent work connected with the Kakeya problem}, Prospects in mathematics (Princeton,
NJ, 1996), 129-162, Amer. Math. Soc., Providence, RI, 1999.


\end{thebibliography}

\end{document}